\documentclass{article}[16pt]
\usepackage{amsfonts,amsmath}
\usepackage[russian]{babel}
\usepackage[cp1251]{inputenc}

\newtheorem{theorem}{Теорема}

\begin{document}

\title{Суммирование рядов Фурье на бесконечномерном торе}

\author{Д.В. Фуфаев \thanks{Механико-математический факультет Московского Государственного Университета имени М.В. Ломоносова,
Россия; e-mail:
fufaevdv@rambler.ru}  }

\date{}

\maketitle

Как известно, даже для функций одной переменной не решена до конца проблема условия сходимости почти всюду рядов Фурье. В то же время для их средних арифметических --- средних Фейера --- ответ известен: достаточным условием сходимости средних Фейера функции $f$ к ней самой почти всюду является суммируемость по Лебегу функции $f$. В случае функций многих переменных и, соответственно, кратных рядов, возникает специфика, связанная с возможностью разных скоростей возрастания индексов суммирования по разным координатам. Отсюда следует возможность рассмотрения разных видов сходимости рядов Фурье и построенных по ним средних. Классически рассматриваемыми случаями являются случай, когда все индексы возрастают одинаково (сходимость по кубам), и случай, когда нет ограничений на скорости возрастания (сходимость по Прингсхейму). Для рядов Фурье вопрос условия сходимости почти всюду также остается открытым, для средних Фейера же ситуация ясна: для сходимости средних Фейера к исходной функции почти всюду по кубам по-прежнему достаточно суммируемости исходной функции по Лебегу, в случае же сходимости по Прингсхейму условием является принадлежность функции классу $L(\ln^+L)^{N-1}$, где $N$ --- размерность пространства (см. \cite{1}, гл.XVII, \S2-3). Более того, эти условия являются неусиляемыми. На самом деле, вместо сходимости по кубам можно рассматривать более слабое условие (но, тем не менее, не ослабляющее описанных результатов) регулярного возрастания индексов, однако, во-первых, на практике этот случай приводит к работе с достаточно громоздкими обозначениями, во-вторых, при обобщении метода суммирования средними арифметическими на другие методы работа с регулярно возрастающими индексами оказывается принципиально невозможной (например, как в случае со средними Марцинкевича). Кроме того, представляют интерес промежуточные случаи регулярности возрастания индексов, между возрастанием по кубам и по Прингсхейму --- тогда получаются соответственно промежуточные условия на функцию для сходимости средних почти всюду.

Естественным переходом от функций многих переменных является рассмотрение функций счетного числа переменных, то есть функций на бесконечномерном торе. Систематическое исследование в этой области провел Йессен (\cite{2}). Для рядов Фурье (точнее --- их бесконечномерного аналога) в этом случае уже имеется окончательное условие для сходимости почти всюду (правда, только по кубам), полученное Холщевниковой Н.Н. (\cite{3}). Анализ на бесконечномерном торе представляет интерес и для гармонического анализа на топологических группах (см. \cite{4}), и для приложений в теории вероятностей и математической физике (см., например, \cite{6}-\cite{9}).

Настоящая работа продолжает исследование анализа на бесконечномерном торе и устанавливает условия для сходимости почти всюду средних Фейера (точнее --- их бесконечномерного аналога) для различных случаев возрастания индексов.

\section{Основные результаты}
\label{subsec2}

Бесконечномерным тором $\mathbb T^{\infty}$ называется декартово произведение счетного числа одномерных торов $\mathbb T$. На $\mathbb T^{\infty}$ определена мера Лебега, как произведение нормированных единицей одномерных мер Лебега на торах-сомножителях. В силу того, что бесконечномерный тор является компактной топологической группой, на нем есть мера Хаара, которая совпадает с упомянутой мерой Лебега. Относительно этой меры вводятся пространства Лебега $L^p(\mathbb T^{\infty})$ и другие.

На $\mathbb T^{\infty}$ рассматривается система функций, называемая системой Йессена:

$$
\theta_{n_1,\dots,n_p}(x)=\prod\limits_{r=1}^{p}e^{2\pi in_rx_r},
$$
где $p\in\mathbb N$, $n_r\in\mathbb Z$, $x=(x_1,\dots,x_p,\dots)\in\mathbb T^{\infty}$ (подробнее об этой системе см. \cite{2}). На самом деле, система Йессена является группой характеров бесконечномерного тора.

Следуя \cite{3}, обозначим через $\mathbb Z^{<\infty}$ множество бесконечномерных целочисленных векторов $n=(n_1,\dots,n_p,\dots)$, у которых лишь конечное число координат отлично от нуля. Сложение векторов и умножение на числа определяются покоординатно.

 Рассмотрим ряд по системе Йессена:

\begin{equation}
\sum\limits_{n\in\mathbb Z^{<\infty}}a_ne^{2\pi inx},
\label{Fuf_eq0}
\end{equation}
где $x\in\mathbb T^{\infty}, nx=\sum\limits_{k=1}^{\infty}n_kx_k,a_n\in\mathbb C$. Для $n\in\mathbb Z^{<\infty}$ будем пользоваться обозначением $a_n=a_{n_1,\dots,n_p}$, если $n_k=0$ при $k>p$.

Прямоугольные частичные суммы такого ряда имеют вид 

\begin{equation}
S_{p,N_1\dots,N_p}(x)=\sum\limits_{n_1=-N_1}^{N_1}\dots\sum\limits_{n_p=-N_p}^{N_p}a_{n_1,\dots,n_p}e^{2\pi i(n_1x_1+\dots+n_px_p)}.
\label{Fuf_eq1}
\end{equation}

Если $N_1=N_2=\dots=N_p$, то имеем кубические частичные суммы.

Ряд  \eqref{Fuf_eq0} называется сходящимся по прямоугольникам в точке $x$ к числу $s$, если для любого $\varepsilon>0$ существует такое $P\in\mathbb N$, что для любого $p>P$ найдется такое $N\in\mathbb N$, что для всех $N_1,\dots,N_p\ge N$ выполнено

$$
|S_{p,N_1\dots,N_p}(x)-s|<\varepsilon.
$$ 

Рассмотрим более сильную сходимость, а именно, скажем, что ряд сходится в усиленном смысле, если ряд $S_{p}=\sum\limits_{n_1,\dots,n_p}a_{n_1,\dots,n_p}e^{2\pi i(n_1x_1+\dots+n_px_p)}$ сходится для каждого $p\in\mathbb N$ и существует конечный предел сумм: $
\lim\limits_{p\to\infty}S_{p}.
$

Подробнее о сходимости рядов на бесконечномерном торе см. \cite{9}.

Через $\mathbb T^m$ будем обозначать $m$-мерный тор, образованный первыми $m$ координатами бесконечномерного, а через $\mathbb T^{m,\infty}$ --- ``хвост'' бесконечномерного тора:

$$
\mathbb T^{m,\infty}=\{ (x_{m+1},\dots,x_{m+p},\dots), x_{m+p}\in\mathbb T, p\in\mathbb N          \}
$$

Пусть $f\in L(\mathbb T^{\infty})$. Положим 

$$
f_m(x)=\int\limits_{\mathbb T^{m,\infty}}f(x)dx_{m+1}\dots dx_{m+p}\dots.
$$

Фактически, эти функции зависят лишь от первых $m$ координат и потому определены на $\mathbb T^m$, но мы будем их рассматривать как функции на $\mathbb T^{\infty}$. По теореме Фубини эти функции существуют для почти всех $(x_1,\dots,x_m)$ и выполнено равенство

$$
\int\limits_{\mathbb T^{\infty}}f(x)dx_{1}\dots dx_{p}\dots =\int\limits_{\mathbb T^{m}}\int\limits_{\mathbb T^{m,\infty}}f(x)dx_{m+1}\dots dx_{m+p}\dots dx_{1}\dots dx_{m}.
$$

Йессен (\cite{2}) доказал, что как функции на бесконечномерном торе $f_m$ сходятся к $f$ почти всюду при $m\to\infty$.

Всякой интегрируемой функции соответствует ее ряд Фурье по системе Йессена:

$$
S(f)=\sum\limits_{n\in\mathbb Z^{<\infty}}c_ne^{2\pi inx}=\sum\limits_{n\in\mathbb Z^{<\infty}}c_n\theta_n(x),\ c_n=\int\limits_{\mathbb T^{\infty}}f(x)e^{-2\pi inx}dx.
$$

Частичные суммы ряда Фурье, то есть суммы вида \eqref{Fuf_eq1} при $a_n=c_n$, можно записать в виде свертки с многомерным ядром Дирихле, то есть в виде

$$
S_{p,N_1\dots,N_p}(x)=\int\limits_{\mathbb T^\infty}\prod\limits_{j=1}^pD_{N_j}(x_j-t_j)f(t)dt,
$$
где $D_l(\cdot)$ --- одномерное ядро Дирихле под номером $l$.
Будем рассматривать средние Фейера на бесконечномерном торе, то есть средние арифметические частичных сумм ряда Фурье по системе Йессена. Их можно записать в виде

$$
\sigma_{p,N_1\dots,N_p}(x)=\int\limits_{\mathbb T^\infty}\prod\limits_{j=1}^pK_{N_j}(x_j-t_j)f(t)dt,
$$
где $K_l(\cdot)=\frac{1}{l+1}\sum\limits_{r=0}^lD_r(\cdot)$ --- одномерное ядро Фейера с номером $l$.

Коэффициенты Фурье функции $f_m$ совпадают с коэффициентами Фурье функции $f$ для индексов $s=(s_1,\dots,s_p,0,\dots,0,\dots)\in\mathbb Z^{<\infty}$ при $p\le m$, и равны нулю, если хотя бы одно из чисел $s_{m+1},s_{m+2},\dots$ не равно нулю. Отсюда следует совпадение частичных сумм ряда Фурье и средних Фейера функции $f_m$  с частичными суммами ряда Фурье и средними Фейера функции $f$.

Из теоремы Фейера (см., например, \cite{1}, гл. XVII, 3.1) следует, что средние Фейера функции $f_m$ сходятся к $f_n$ почти всюду на $\mathbb T^m$ по кубам. Отсюда и из теоремы Йессена следует 

\begin{theorem} Пусть $ f \in L(\mathbb T^\infty)$. Тогда средние Фейера сходятся в усиленном смысле к $f$ по кубам почти всюду на $\mathbb T^\infty$.

\end{theorem}

{\textbf {Замечание.}} Данную теорему можно обобщить, если рассматривать сходимость не по кубическим индексам, а сходимость при регулярном возрастании индексов, так как утверждение теоремы Фейера остается верным для этого случая (регулярность здесь понимается в том же смысле, что и ограниченность, см.\cite{1}, гл.XVII, \S 3).


Для пространства с мерой $(X,\mu)$ через $L(\ln^{+}L)^D  (X), D\in \mathbb N,$ обозначается пространство измеримых комплекснозначных функций $f$, удовлетворяющих условию

$$
\int\limits_{X}|f(x)|\ln^D(|f(x)|+1)d\mu(x)<\infty.
$$

Чтобы получить результаты о сходимости по Прингсхейму, нам понадобится слещующее утверждение.

{\textbf {Лемма.}} Пусть $ f \in L(\ln^+L)^d(\mathbb T^\infty)$, $d\in\mathbb N$. Тогда для любого $m\in\mathbb N$ имеем $f_m\in L(\ln^+L)^d(\mathbb T^m)$. Кроме того, если $ f \in \bigcap\limits_{r=1}^{\infty}L(\ln^+L)^r(\mathbb T^\infty)$, то для любых $m,d\in\mathbb N$ имеем  $f_m\in L(\ln^+L)^d(\mathbb T^m)$.

 {\textbf {Доказательство.}}
Напомним следующее неравенство Йенсена для интеграла (см. \cite{10}, гл. 3, теорема 3.3): для действительнозначной функции $g$, интегрируемой на пространстве с мерой $(X,\nu)$, $\nu(X)=1$, и выпуклой вниз на области значений функции $g$ действительнозначной функции $\phi$ справедливо неравенство

$$
\phi \left(\int\limits_{X}g(x)d\nu(x) \right) \le \int\limits_{X}\phi(g(x))d\nu(x).
$$

Применим его к случаю $X=\mathbb T^{m,\infty}$, $\phi(x)=x\ln(x+1)$, $g(x)=|f(x)|$.

$$
|f_m(x)| (\ln(|f_m(x)|+1))^d\le\int\limits_{\mathbb T^{m,\infty}}|f(x)|dx_{m+1}\dots\times \left( \ln\left(\int\limits_{\mathbb T^{m,\infty}}|f(x)|dx_{m+1}\dots+1\right)    \right)^d\le
$$

$$
\le\int\limits_{\mathbb T^{m,\infty}}|f(x)|(\ln(|f(x)|+1))^ddx_{m+1}\dots.
$$

Здесь также использовали возрастание функции $\phi$. Проинтегрируем последнее неравенство по $\mathbb T^m$:

$$
\int\limits_{\mathbb T^m}|f_m(x)| (\ln(|f_m(x)|+1))^ddx_{1}\dots dx_m\le\int\limits_{\mathbb T^{\infty}}|f(x)|(\ln(|f(x)|+1))^ddx_{1}\dots,
$$
откуда $f_m\in L(\ln^+L)^d(\mathbb T^m)$.

\begin{theorem} Пусть $ f \in \bigcap\limits_{r=1}^{\infty}L(\ln^+L)^r(\mathbb T^\infty)$. Тогда средние Фейера сходятся в усиленном смысле к $f$ по Прингсхейму почти всюду на $\mathbb T^\infty$.

\end{theorem}

 {\textbf {Доказательство.}}

Зафиксируем $m\in\mathbb N$. Из леммы следует, что $f_m\in L(\ln^+L)^{m-1}(\mathbb T^m)$. Тогда по теореме Йессена-Марцинкевича-Зигмунда (\cite{1}, гл. XVII, 2.14) средние Фейера функции $f_m$ сходятся к $f_m$ почти всюду, откуда, по теореме Йессена, следует утверждение теоремы.

{\textbf {Замечание.}} В силу неусиляемости теоремы Йессена-Марцинкевича-Зигмунда в конечномерном случае, в каждом классе $ L(\ln^+L)^{d}(\mathbb T^\infty)$ найдется функция, средние Фейера которой расходятся всюду по Прингсхейму.

Рассмотрим теперь промежуточный случай регулярности возрастания индексов, т.е. $D$-регулярное возрастание (см. \cite{11}). Здесь это можно сформулировать так: пусть бесконечномерный тор представлен в виде произведения некоторого числа торов: $\mathbb T^\infty=\prod\limits_{k=1}^D\mathbb T^{M_k}$, где $M_k\in\mathbb N\cup\{\infty\}, D\in\mathbb N$. Обозначим через $P^k$ проектор на соответствующее подпространство координат. Тогда индекс $n$ возрастает $D$-регулярно, если для каждого $k\in\{1,\dots,D\}$ регулярно возрастает индекс $P^k(n)\in\mathbb Z^{M_k}$.

\begin{theorem} Пусть $ f \in L(\ln^+L)^{D-1}(\mathbb T^\infty)$. Тогда средние Фейера сходятся в усиленном смысле к $f$ $D$-регулярно почти всюду на $\mathbb T^\infty$.

\end{theorem}

 {\textbf {Доказательство.}}
Применяя лемму, получаем, что $f_m\in L(\ln^+L)^{D-1}(\mathbb T^m)$ для любого $m\in\mathbb N$. Тогда по доказанному в \cite{11} средние Фейера функции $f_m$ сходятся к $f_m$ $D$-регулярно почти всюду. Применяя теорему Йессена, аналогично предыдущим случаям получаем утверждение теоремы.

{\textbf {Замечание.}} В этой теореме можно допускать случай $D=\infty$ (требование на функцию тогда будет: $ f \in \bigcap\limits_{r=1}^{\infty}L(\ln^+L)^r(\mathbb T^\infty)$), который доказывается аналогично предыдущим теоремам. Этот случай можно считать некоторым обобщением теоремы 2, ведь $\infty$-регулярность требует тех же условий, что и сходимость по Прингсхейму, однако здесь допускается случай, когда ``кубические'' индексы не одномерны, а многомерны, и даже бесконечномерны, причем здесь допускается даже бесконечное число бесконечномерных торов в произведении.

{\textbf {Замечание.}} Вместо средних Фейера можно рассматривать средние рядов Фурье более общего вида, которые будут гарантировать сходимость к разлагаемой функции почти всюду (см., например, \cite{12}-\cite{14}).

\section{Абстрактные пространства с мерой}
\label{subsec3}

Полученные результаты можно обобщить на ситуацию более общего характера, а именно: пусть $(X^k,\mu^k)$, $k\in\mathbb N$ --- пространства с мерой, причем $\mu^k(X^k)=1$.

Для дальнейшего необходимо напомнить определение направленности (см., например, \cite{15}, т.2, стр.12). Непустое множество $A$ называется направленным, если на нем задан частичный порядок, удовлетворяющий следующему условию: для любых $q,s\in A$ найдется элемент $n\in A$ такой, что $q\le n$ и $s\le n$. Направленностью в множестве $X$ называется набор элементов $\{x_n\}_{n\in A}$, индексируемых элементами направленного множества. Направленность $\{x_n\}_{n\in A}$ в топологическом пространстве $X$ сходится к элементу $x$, если для любого непустого открытого множества $U$, содержащего $x$, найдется такой элемент $n_0\in A$, что $x_n\in U$ для всех $n\ge n_0, n\in A$. Понятным образом определяется сходимость числовых направленностей, а также поточечная сходимость и сходимость почти всюду направленностей числовых функций.

Пусть теперь для каждого $k\in\mathbb N$ задана направленность операторов $\{T^k_{n^k}\}_{n^k\in A^k}$, переводящая $L^1(X^k,\mu^k)$ в себя, причем в каждом направленном множестве существует ``нулевой'' элемент $0^k\in A^k$, наименьший по порядку, такой, что соответствующий оператор направленности есть тривиальный оператор интегрирования: $T^k_{0^k}f=I_{X^k}(f)=\int\limits_{X^k}f(x)d\mu^k(x)$ для любой функции $f\in L^1(X^k,\mu^k)$. От разрозненных пространств функций и операторов на них нам надо перейти к функциям, заданным на произведении пространств. Для этого напомним понятие тензорного произведения (\cite{16}, 10.42.1, 10.42.2).

 Пусть $f^1\in L^1(X^1,\mu^1)$, $f^2\in L^1(X^2,\mu^2)$, тогда их тензорным произведением называется функция $f^1\otimes f^2\in L^1(X^1\times X^2, \mu^1\otimes \mu^2)$, определяемая по формуле $f^1\otimes f^2(x_1,x_2)=f^1(x_1)f^2(x_2)$, где $\mu^1\otimes \mu^2$ --- это произведение мер $\mu^1$ и $\mu^2$ (см., например, \cite{15}, т.1, стр.223). Такие функции называются элементарными тензорами. Алгебраическим тензорным произведением $L^1(X^1,\mu^1)$ и $L^1(X^2,\mu^2)$ называется пространство линейных комбинаций элементарных тензоров и обозначается $L^1(X^1,\mu^1) \mathbin{\otimes} L^1( X^2,\mu^2)$. Проективным тензорным произведением, обозначаемым $L^1(X^1,\mu^1) \mathbin{\hat{\otimes}} L^1( X^2,\mu^2)$, называется пополнение алгебраического по 
так называемой проективной норме (\cite{17}, 0.3.28), но нам она не нужна, а нужен тот факт, что оно
изометрически изоморфно $L^1(X^1\times X^2,\mu^1\otimes \mu^2)$ (\cite{17}, следствие 0.3.36).

Далее, если непрерывные линейные операторы $T^i$ переводят $L^1(X^i,\mu^i)$ в себя, $i=1,2$, то существует единственный непрерывный линейный оператор $T^1\mathbin{\hat{\otimes}}T^2$, переводящий $L^1(X^1,\mu^1) \mathbin{\hat{\otimes}} L^1( X^2,\mu^2)$ в себя, такой, что на элементарных тензорах он действует по формуле: $T^1\mathbin{\hat{\otimes}}T^2[f\cdot g](x,y)=T^1f(x)\cdot T^2g(y)$ (\cite{17}, теорема 0.3.40).

Вообще говоря, этого недостаточно для того, чтобы определить тензорное произведение бесконечного числа пространств и, тем более, операторов. Но нам достаточно будет считать, что тензорное произведение бесконечного числа тривиальных операторов --- тривиальный оператор, то есть, что $\hat{\bigotimes\limits_{k\ge j}}T^k_{0^k}=I_{\prod\limits_{k\ge j}X^k}$ для любого $j\in\mathbb N$ (о произведении бесконечного числа пространств с мерой см. \cite{15}, т.1, 3.5). Тогда мы можем определить бесконечное произведение операторов, лишь конечное число из которых отлично от тривиального: $T^1_{n^1}\mathbin{\hat{\otimes}}\dots\mathbin{\hat{\otimes}}T^p_{n^p}\mathbin{\hat{\otimes}}T^{p+1}_{0^{p+1}}\mathbin{\hat{\otimes}}\dots=T^1_{n^1}\mathbin{\hat{\otimes}}\dots\mathbin{\hat{\otimes}}T^p_{n^p}\mathbin{\hat{\otimes}}I_{\prod\limits_{k\ge p+1}X^k}$ для любого $p\in\mathbb N$. Последнее выражение определено корректно.

Для двух направленностей операторов $\{T^1_{n^1}\}_{n^1\in A^1}$ и $\{T^2_{n^2}\}_{n^2\in A^2}$ их тензорным произведением назовем направленность
$\{T^1_{n^1} \mathbin{\hat{\otimes}}T^2_{n^2} \}_{\substack {   n\in A }} $, где $ n=(n^1,n^2), A=A^1\times A^2$, причем $(n^1,n^2)>(m^1,m^2)$ тогда и только тогда, когда $n^1>m^1$ и $n^2>m^2$.
Тензорное произведение любого конечного числа направленностей определяется аналогичным образом. Тензорное произведение бесконечного числа направленностей операторов определяется как направленность $\{T^1_{n^1}\mathbin{\hat{\otimes}}\dots\mathbin{\hat{\otimes}}T^k_{n^k}\mathbin{\hat{\otimes}}\dots \}_{\substack {    n\in A }}$, где $ n=(n^1,\dots,n^k,\dots)$, $A\subset\prod\limits_{j=1}^{\infty}A^j$, направленное множество $A$ состоит из последовательностей, у которых лишь конечное число элементов отлично от ``нулевого'', то есть от $0^k, k\in\mathbb N$. 

Согласно \cite{15}, т.2, 10.2.4, последовательность $\text{id}_{L^1\left(\prod\limits_{k=1}^p X^k, \bigotimes\limits_{k=1}^p\mu^k\right)}\mathbin{\hat{\otimes}}I_{\prod\limits_{k> p+1}X^k}f$ сходится $ \bigotimes\limits_{k=1}^\infty\mu^k$-почти всюду к $f$ при $p\to\infty$ для любой $f\in L^1\left(\prod\limits_{k=1}^\infty X^k, \bigotimes\limits_{k=1}^\infty\mu^k\right)$. Теперь теоремы 1-3 обобщаются следующим образом (доказательства аналогичны): 

\begin{theorem} {Пусть $(X^k,\mu^k)$, $k\in\mathbb N$ --- пространства с мерой, $\mu_k(X^k)=1$, для каждого $k\in\mathbb N$ задана направленность операторов $\{T^k_{n^k}\}_{n^k\in A^k}$, переводящая $L^1(X^k,\mu^k)$ в себя, каждая из которых имеет ``нулевой'' элемент, тогда

а) если для любого $p\in\mathbb N$ и любой функции $f\in L^1\left(\prod\limits_{k=1}^p X^k, \bigotimes\limits_{k=1}^p\mu^k\right)$ направленность  $\hat{\bigotimes\limits_{1\le k\le p}}T^k_{n^k}f$ сходится $\bigotimes\limits_{k=1}^p\mu^k$-почти всюду к $f$, то для любой функции $f\in L^1\left(\prod\limits_{k=1}^\infty X^k, \bigotimes\limits_{k=1}^\infty\mu^k\right)$ направленность  $\hat{\bigotimes\limits_{1\le k\le \infty}}T^k_{n^k}f$ сходится $\bigotimes\limits_{k=1}^\infty\mu^k$-почти всюду к $f$.

б) если для любого $p\in\mathbb N$ и любой функции $f\in  L(\ln^+L)^{p-1}\left(\prod\limits_{k=1}^p X^k, \bigotimes\limits_{k=1}^p\mu^k\right)$ направленность  $\hat{\bigotimes\limits_{1\le k\le p}}T^k_{n^k}f$ сходится $\bigotimes\limits_{k=1}^p\mu^k$-почти всюду к $f$, то для любой функции $f\in \bigcap\limits_{r=1}^{\infty}L(\ln^+L)^r\left(\prod\limits_{k=1}^\infty X^k, \bigotimes\limits_{k=1}^\infty\mu^k\right)$ направленность  $\hat{\bigotimes\limits_{1\le k\le \infty}}T^k_{n^k}f$ сходится к $f$ \newline $\bigotimes\limits_{k=1}^\infty\mu^k$-почти всюду.}
\end{theorem}

В качестве частного случая можно брать не более чем счетные направленности интегральных операторов с мажорантой слабого типа (1,1), см. \cite{11}.

\end{document}